# MODERATE DEVIATIONS FOR DIFFUSIONS WITH BROWNIAN POTENTIALS

By Yueyun Hu and Zhan Shi

*Université Paris VI*

We present precise moderate deviation probabilities, in both quenched and annealed settings, for a recurrent diffusion process with a Brownian potential. Our method relies on fine tools in stochastic calculus, including Kotani's lemma and Lamperti's representation for exponential functionals. In particular, our result for quenched moderate deviations is in agreement with a recent theorem of Comets and Popov [*Probab. Theory Related Fields* **126** (2003) 571–609] who studied the corresponding problem for Sinai's random walk in random environment.

**1. Introduction.** Let $W := (W(x), x \in \mathbb{R})$ be a one-dimensional Brownian motion defined on $\mathbb{R}$ with $W(0) = 0$. Let $(\beta(t), t \geq 0)$ be another one-dimensional Brownian motion independent of $W(\cdot)$. Following [2] and [20], we consider the equation $X(0) = 0$,

(1.1) $$dX(t) = d\beta(t) - \tfrac{1}{2}W'(X(t))\,dt, \qquad t \geq 0.$$

The solution $X$ of (1.1) is called a diffusion with random potential $W$. The rigorous meaning of (1.1) can be given in terms of the infinitesimal generator: Conditioning on each realization $\{W(x), x \in \mathbb{R}\}$, the process $X$ is a real-valued diffusion with generator

$$\frac{1}{2}e^{W(x)}\frac{d}{dx}\left(e^{-W(x)}\frac{d}{dx}\right).$$

Another representation of $X$ by time change is given in Section 4.

The process $X$ has been used in modeling some random phenomena in physics [16]. It is also related to random walk in random environment [18, 7, 26]. See [24] and [21] for recent surveys.

We denote by $P$ and $E$ the probability and the expectation with respect to the potential $W$, and by $\mathbf{P}_\omega$ and $\mathbf{E}_\omega$ the quenched probability and the quenched expectation ("quenched" means the conditioning with respect to the potential $W$). The total (or annealed) probability is $\mathbb{P} \stackrel{\text{def}}{=} P(d\omega) \otimes \mathbf{P}_\omega$.

The typical long-time behavior of $X(t)$ is described by a result of Brox [2], which is the continuous-time analogue of Sinai's [22] well-known theorem for recurrent random walk in random environment: under the total probability $\mathbb{P}$,

$$\frac{X(t)}{\log^2 t} \stackrel{(d)}{\to} b(1), \qquad t \to \infty,$$

where $\stackrel{(d)}{\to}$ denotes convergence in distribution, and $b(1)$ is a nondegenerate random variable whose distribution is explicitly known.

It is interesting to study the deviation probabilities

(1.2)     $\mathbf{P}_\omega\{X(t) > v\}$  and   $\mathbb{P}\{X(t) > v\}, \qquad t, v \to \infty, v \gg \log^2 t,$

where $v \gg \log^2 t$ means $v/\log^2 t \to \infty$. In the sequel, we also write $x \ll y$ or $x = o(y)$ to denote $y \gg x$. When $v/t$ converges to a positive constant, this is a large deviation problem, and is solved by Taleb [23] (who actually studies the problem for all drifted Brownian potentials). In particular, it is shown that in this case both probabilities in (1.2) have exponential decays.

We focus on moderate deviation probabilities, that is, when $(v, t)$ is such that $\log^2 t \ll v \ll t$.

Our first result, which concerns the quenched setting, is in agreement with Theorem 1.2 of [3] for random walk in random environment. This was, indeed, the original motivation of the present work.

THEOREM 1.1. *We have*

(1.3)     $$\frac{2\log(t/v)}{v} \log \mathbf{P}_\omega\{X(t) > v\} \to -1, \qquad P\text{-a.s.},$$

*whenever $v, t \to \infty$ such that $v \gg (\log^2 t) \log \log \log t$ and $\log \log t = o(\log(t/v))$. The same result holds for $\sup_{0 \le s \le t} X(s)$ instead of $X(t)$.*

Loosely speaking, Theorem 1.1 says that in a typical potential $W$, $\mathbf{P}_\omega\{X(t) > v\}$ behaves like $\exp[-(1+o(1))\frac{v}{2\log(t/v)}]$. However, if we take the average over all the realizations of $W$ (i.e., in the annealed setting), the deviation probability will become considerably larger. This is confirmed in our second result stated as follows.

THEOREM 1.2. *We have*

(1.4)     $$\frac{\log^2 t}{v} \log \mathbb{P}\{X(t) > v\} \to -\frac{\pi^2}{8},$$

*whenever* $v, t \to \infty$ *such that* $v \gg \log^2 t$ *and* $\log v = o(\log t)$. *The same result holds for* $\sup_{0 \le s \le t} X(s)$ *instead of* $X(t)$.

When $\log^2 t \le v \le (\log^2 t)(\log \log t)^{1/2}$, the convergence (1.4) has already been obtained in [9] by means of the Laplace method. This method, however, fails when $v$ goes to infinity too quickly, for example, if $v \gg \log^3 t$.

We say a few words about the proofs of Theorems 1.1 and 1.2. Although the methods adopted in the three parts (Theorem 1.1, upper bound and lower bound in Theorem 1.2) find all their roots in stochastic calculus, they rely on completely different ingredients.

In the proof of Theorem 1.1, we exploit Kotani's lemma as well as other fine tools in the theory of one-dimensional diffusion.

The proof of the upper bound in Theorem 1.2 relies on Lamperti's representation for exponential functionals and on Warren and Yor's [25] skew-product theorem for Bessel processes. The proof of the lower bound, on the other hand, is based on a bare-hand analysis of pseudo-valleys where the diffusion $X$ spends much time.

The rest of the paper is organized as follows. Section 2 is devoted to some preliminary results for local times of Brownian motion. In Section 3 we introduce Kotani's lemma and prove Theorem 1.1. The main result in Section 4, Theorem 4.1, is a joint arcsine type law for the occupation times of $X$, which may be of independent interest. This result will be used to prove Theorem 1.2 in Section 5.

Throughout this paper, we write $f^{-1}$ for the inverse of any continuous and (strictly) increasing function $f$. Unless stated otherwise, for any continuous process $\xi$, we denote by $T_\xi(x) = \inf\{t \ge 0 : \xi(t) = x\}$, $x \in \mathbb{R}$, the first hitting time of $\xi$ at $x$.

**2. Preliminaries on local times.** In this section we collect a few preliminary results for the local times of Brownian motion. These results will be of use in the rest of the paper.

Let $B$ be a one-dimensional Brownian motion starting from 0. Let $(L(t, x), t \ge 0, x \in \mathbb{R})$ be the family of the local times of $B$, that is, for any Borel function $f : \mathbb{R} \to \mathbb{R}_+$, $\int_0^t f(B(s))\,ds = \int_{-\infty}^\infty f(x) L(t, x)\,dx$. We define

$$\tau(r) \stackrel{\text{def}}{=} \inf\{t > 0 : L(t, 0) > r\}, \qquad r \ge 0, \tag{2.1}$$

$$\sigma(x) \stackrel{\text{def}}{=} \inf\{t > 0 : B(t) > x\}, \qquad x \ge 0. \tag{2.2}$$

Denote by BES($\delta$) [resp. BESQ($\delta$)] the Bessel process (resp. the squared Bessel process) of dimension $\delta$. We recall that a $\delta$-dimensional squared Bessel process has generator of form $2x \frac{d^2}{dx^2} + \delta \frac{d}{dx}$. When $\delta$ is an integer, a Bessel process can be realized as the Euclidean norm of an $\mathbb{R}^\delta$-valued Brownian motion. We refer to [19], Chapter XI, for a detailed account of general properties of Bessel and squared Bessel processes, together with the proof of the following result.

FACT 2.1 (First Ray–Knight theorem). *Fix $a > 0$. The process $\{L(\sigma(a), a - x), x \geq 0\}$ is an inhomogeneous strong Markov process starting from 0, which is a BESQ(2) on $[0, a]$ and a BESQ(0) on $(a, \infty)$.*

The rest of the section is devoted to a few preliminary results for local times of Brownian motion.

LEMMA 2.2. *For $b > a \geq 0$ and $v > 0$,*

(2.3)
$$\frac{2}{\pi} \exp\left(-\frac{\pi^2}{8} \frac{v}{(b-a)^2}\right) \leq \mathbb{P}\left(\int_a^b L(\sigma(b), x)\, dx > v\right)$$
$$\leq \frac{4}{\pi} \exp\left(-\frac{\pi^2}{8} \frac{v}{(b-a)^2}\right).$$

PROOF. By the strong Markov property, $\int_a^b L(\sigma(b), x)\, dx$ is distributed as $\int_0^{b-a} L(\sigma(b-a), x)\, dx$. According to Fact 2.1 and [19], Corollary XI.1.8, for any $\lambda > 0$, $\mathbb{E}\exp\{-\frac{\lambda^2}{2} \int_0^{b-a} L(\sigma(b-a), x)\, dx\} = 1/\cosh(\lambda(b-a))$, which is also the Laplace transform at $\frac{\lambda^2}{2}$ of $T_{|B|}(b-a)$ (the first hitting time of $b - a$ by $|B|$). Thus, $\int_a^b L(\sigma(b), x)\, dx \stackrel{\text{law}}{=} T_{|B|}(b-a)$. The lemma now follows from the exact distribution of Brownian exit time ([6], page 342).  □

LEMMA 2.3. *Let $b > a > 0$ and $\kappa > 0$, we have*
$$\int_0^{\sigma(a)} (b - B(s))^{(1/\kappa)-2}\, ds \stackrel{\text{law}}{=} \Upsilon_{2-2\kappa}(2\kappa b^{1/(2\kappa)} \rightsquigarrow 2\kappa(b-a)^{1/(2\kappa)}),$$

*where $\Upsilon_{2-2\kappa}(x \rightsquigarrow y)$ means the first hitting time of $y$ by a BES$(2 - 2\kappa)$ starting from $x$.*

PROOF. Write
$$\int_0^{\sigma(b)} (b - B(s))^{(1/\kappa)-2}\, ds$$
$$= \int_0^{\sigma(a)} (b - B(s))^{(1/\kappa)-2}\, ds + \int_{\sigma(a)}^{\sigma(b)} (b - B(s))^{(1/\kappa)-2}\, ds.$$

By the strong Markov property, the integrals on the right-hand side are independent random variables. Moreover, the second integral is distributed as $(b-a)^{1/\kappa} \int_0^{\sigma(1)} (1 - B(s))^{(1/\kappa)-2}\, ds$.

On the other hand, according to Getoor and Sharpe ([8], Proposition 5.14(a)), for any $\lambda > 0$,

$$\mathbb{E}\exp\left(-\frac{\lambda^2}{2}\int_0^{\sigma(1)}(1-B(s))^{(1/\kappa)-2}\,ds\right) = \frac{2}{\Gamma(\kappa)}(\lambda\kappa)^\kappa K_\kappa(2\kappa\lambda),$$

where $K_\kappa$ denotes the modified Bessel function of index $\kappa$.

Assembling these pieces yields that

$$\mathbb{E}\exp\left(-\frac{\lambda^2}{2}\int_0^{\sigma(a)}(b-B(s))^{(1/\kappa)-2}\,ds\right) = \left(\frac{b}{b-a}\right)^{1/2}\frac{K_\kappa(2\kappa\lambda b^{1/(2\kappa)})}{K_\kappa(2\kappa\lambda(b-a)^{1/(2\kappa)})}.$$

According to Kent [13], the expression on the right-hand side is exactly the Laplace transform at $\lambda^2/2$ of $\Upsilon_{2-2\kappa}(2\kappa b^{1/(2\kappa)} \rightsquigarrow 2\kappa(b-a)^{1/(2\kappa)})$. □

LEMMA 2.4. *Almost surely,*

$$\frac{1}{r}\sup_{|x|\leq u}(L(\tau(r),x)-r)\to 0,$$

*whenever $u\to\infty$ and $r \gg u\log\log u$.*

PROOF. By symmetry, we only need to treat the case $0 \leq x \leq u$. It is proved by Csáki and Földes ([4], Lemma 2.1) that for any $\varepsilon > 0$,

$$\lim_{r\to\infty}\frac{1}{r}\sup_{0\leq x\leq r/(\log\log r)^{1+\varepsilon}}(L(\tau(r),x)-r)=0 \quad \text{a.s.}$$

Thus, we only have to deal with the case $r \gg u\log\log u$ and $r \leq u\log u$.

Let $\varepsilon > 0$. We shall prove that almost surely for all $u,r\to\infty$ such that $r \gg u\log\log u$ and $r \leq u\log u$,

(2.4) $$\frac{1}{r}\sup_{0\leq x\leq u}(L(\tau(r),x)-r)\leq \varepsilon.$$

To this end, we consider the events

$$A_{k,j}\stackrel{\text{def}}{=}\left\{\exists\,(u,r)\in[u_k,u_{k+1}]\times[r_j,r_{j+1}]:\frac{512}{\varepsilon^2}u\log\log u\leq r\leq u\log u,\right.$$

$$\left.\sup_{0\leq x\leq u}(L(\tau(r),x)-r)>\varepsilon r\right\},$$

with $u_k\stackrel{\text{def}}{=}2^k$, $r_j\stackrel{\text{def}}{=}2^j$ and $k,j\geq 100$. The desired conclusion (2.4) will follow from the Borel–Cantelli lemma once we show that

$$\sum_{j,k\geq 100}\mathbb{P}(A_{k,j})<\infty.$$

To prove $\sum_{j,k\geq 100} \mathbb{P}(A_{k,j}) < \infty$, we recall a result of Bass and Griffin ([1], Lemma 3.4) saying that there exists a constant $c > 0$ such that for all $0 < h, x < 1$,

$$\mathbb{P}\left\{\sup_{1\leq r\leq 2}\sup_{0\leq z\leq h}(L(\tau(r),z)-r) > x\right\} \leq c\frac{x}{h}\exp\left(-\frac{x^2}{32h}\right).$$

Therefore, by the monotonicity and the scaling property,

$$\mathbb{P}(A_{k,j}) \leq \mathbb{P}\left(\sup_{0\leq x\leq u_{k+1}}\sup_{r_j\leq r\leq r_{j+1}}(L(\tau(r),x)-r) > \varepsilon r_j\right)$$

$$= \mathbb{P}\left(\sup_{0\leq x\leq u_{k+1}/r_j}\sup_{1\leq r\leq 2}(L(\tau(r),x)-r) > \varepsilon\right)$$

$$\leq c\frac{\varepsilon r_j}{u_{k+1}}\exp\left(-\frac{\varepsilon^2 r_j}{32u_{k+1}}\right)$$

$$\leq c\frac{\varepsilon r_j}{u_{k+1}}(k\log 2)^{-4}.$$

To ensure $\mathbb{P}(A_{k,j}) > 0$, necessarily $r_j \leq u_{k+1}\log u_{k+1}$, this implies that $j \leq 2k$. Hence,

$$\sum_{j,k\geq 100}\mathbb{P}(A_{k,j}) \leq \sum_{k\geq 100} c\varepsilon 2k(k\log 2)^{-4}(k+1)\log 2 < \infty,$$

as desired. $\square$

LEMMA 2.5. *Let $r > 0$. We have*

$$(2.5)\quad \mathbb{E}\exp\left(-\lambda\int_0^\infty e^{-s}L(\sigma(r),-s)\,ds\right) = \frac{1}{1+r\sqrt{2\lambda}I_1(\sqrt{8\lambda})}, \qquad \lambda > 0,$$

*where $I_1(\cdot)$ denotes the modified Bessel function of index 1. Consequently, there exists some constant $c > 0$ such that for all $r > c$, $\lambda > 0$ and $0 < a \leq r$, we have*

$$(2.6)\quad \mathbb{P}\left(\int_{-\infty}^a e^{-|x|}L(\sigma(r),x)\,dx > \lambda\right) \leq 3\exp\left(-\frac{\lambda}{8r}\right) + 2\exp\left(-\frac{\lambda}{4ar}\right).$$

PROOF. By Fact 2.1, $s \mapsto L(\sigma(r),-s)$ for $s \geq 0$ is a BESQ(0) starting from $L(\sigma(r),0) \stackrel{\text{law}}{=} 2r\mathbf{e}$, where $\mathbf{e}$ is an exponential variable with mean 1. Let $U_a$ be a BESQ(0) starting from $a > 0$. The Laplace transform of $\int_0^\infty U_a(s)e^{-s}\,ds$ is given by the solution of a Sturm–Liouville equation, see [17]: for all $\lambda > 0$,

$$(2.7)\quad \mathbb{E}\exp\left(-\lambda\int_0^\infty e^{-s}U_a(s)\,ds\right) = \exp\left(\frac{a}{2}\psi'_+(0)\right),$$

where $\psi'_+(0)$ is the right-derivative of the convex function $\psi$ at 0, and $\psi$ is the unique solution, decreasing, nonnegative, of the Sturm–Liouville equation: $\psi(0) = 1$,

$$\tfrac{1}{2}\psi''(x) = \lambda e^{-x}\psi(x), \qquad x > 0.$$

Elementary computations ([17], page 435) show that

$$\psi(x) = I_0(\sqrt{8\lambda}e^{-x/2}), \qquad x \geq 0,$$

which implies that $\psi'_+(0) = -\sqrt{2\lambda}I_1(\sqrt{8\lambda})$, where $I_1$ denotes the modified Bessel function of index 1. Plugging this into (2.7) and integrating with respect to $\mathbf{e}$ give the Laplace transform (2.5). By analytic continuation, for all sufficiently small $\lambda > 0$ (how small depending on $r$),

$$\mathbb{E}\exp\left(\lambda \int_0^\infty e^{-s}L(\sigma(r),-s)\,ds\right) = \frac{1}{1 - r\sqrt{2\lambda}J_1(\sqrt{8\lambda})},$$

where $J_1(\cdot)$ is the Bessel function of index 1. Since $J_1(x) \sim x/2$ when $x \to 0$, there exists some large $c > 0$ such that for all $r > c$, we have

$$\mathbb{E}\exp\left(\frac{1}{4r}\int_0^\infty e^{-s}L(\sigma(r),-s)\,ds\right) = \frac{1}{1 - \sqrt{r/2}J_1(\sqrt{2/r})} < 3.$$

This implies that for all $r > c$ and $\lambda > 0$, we have

$$(2.8) \qquad \mathbb{P}\left(\int_0^\infty e^{-s}L(\sigma(r),-s)\,ds > \lambda\right) \leq 3\exp\left(-\frac{\lambda}{4r}\right).$$

On the other hand, $\sup_{0 \leq x \leq a} L(\sigma(r),x)$ is the maximum of a BESQ(2) over $[r-a,r]$. It follows from reflection principle that

$$\mathbb{P}\left(\int_0^a L(\sigma(r),x)e^{-|x|}\,dx \geq \lambda\right) \leq \mathbb{P}\left(\sup_{0 \leq x \leq a} L(\sigma(r),x) \geq \frac{\lambda}{a}\right)$$

$$\leq 2\mathbb{P}\left(L(\sigma(r),0) \geq \frac{\lambda}{a}\right)$$

$$= 2\exp\left(-\frac{\lambda}{2ar}\right),$$

since $L(\sigma(r),0) \stackrel{\text{law}}{=} 2r\mathbf{e}$. This, together with (2.8), yields (2.6) by triangular inequality. □

**3. The quenched case: proof of Theorem 1.1.** This section is devoted to the proof of Theorem 1.1, by means of the so-called Kotani's lemma. Let $X$ be the diffusion process in a Brownian potential $W$ as in (1.1). We define

$$H(v) \stackrel{\text{def}}{=} \inf\{t > 0 : X(t) > v\}, \qquad v \geq 0.$$

Kotani's lemma (see [12]) gives the Laplace transform of $H(v)$ under the quenched probability $\mathbf{P}_\omega$.

FACT 3.1 (Kotani's lemma). For $\lambda > 0$ and $v \geq 0$, we have

$$\mathbf{E}_\omega(e^{-\lambda H(v)}) = \exp\left(-2\lambda \int_0^v Z(s)\,ds\right),$$

where $Z(\cdot) = Z_\lambda(\cdot)$ is the unique stationary and positive solution of the equation

$$dZ(t) = Z(t)\,dW(t) + (1 + \tfrac{1}{2}Z(t) - 2\lambda Z^2(t))\,dt.$$

Before starting the proof of Theorem 1.1, we study the almost sure behaviors of $\int_0^v Z(s)\,ds$ as $v \to \infty$ and $\lambda \to 0$.

LEMMA 3.2. *We have*

$$\int_0^v Z(s)\,ds = (1 + o(1))\frac{v}{4\lambda \log(1/\lambda)} \qquad a.s., \tag{3.1}$$

*whenever* $\lambda \to 0$ *and* $v \to \infty$ *such that* $v \gg \log^2(1/\lambda) \log\log\log(1/\lambda)$.

PROOF. Without loss of generality, we assume $Z(0) = 1$. Let

$$S(x) \stackrel{\text{def}}{=} \int_1^x e^{(2/y) + 4\lambda y}\,\frac{dy}{y}, \qquad x > 0,$$

which is a scale function of the diffusion $Z$. By Feller's time change representation, there exists a standard one-dimensional Brownian motion $B$, such that

$$Z(t) = S^{-1}(B(\Phi^{-1}(t))), \qquad t \geq 0, \tag{3.2}$$

where

$$\Phi(t) \stackrel{\text{def}}{=} \int_0^t \frac{ds}{h^2(B(s))}, \qquad t \geq 0,$$

$$h(x) \stackrel{\text{def}}{=} S'(S^{-1}(x))S^{-1}(x) = \exp\left(\frac{2}{S^{-1}(x)} + 4\lambda S^{-1}(x)\right), \qquad x \in \mathbb{R},$$

$S^{-1}$ and $\Phi^{-1}$ being the inverse functions of $S$ and $\Phi$, respectively.

Let $L(\cdot, \cdot)$ denote, as before, the local time of $B$, and let $\tau(\cdot)$ be the inverse local time at 0 as in (2.1). We define, for any fixed $\gamma \in [0, 1]$,

$$D_\gamma(r) \stackrel{\text{def}}{=} \int_{-\infty}^\infty (S^{-1}(x))^\gamma \tag{3.3}$$
$$\times \exp\left(-\frac{4}{S^{-1}(x)} - 8\lambda S^{-1}(x)\right) L(\tau(r), x)\,dx, \qquad r \geq 0.$$

We claim that

$$(3.4) \quad D_\gamma(r) = \begin{cases} (1+o(1))r\log(1/\lambda), & \text{if } \gamma = 0, \\ (1+o(1))\Gamma(\gamma)\dfrac{r}{(4\lambda)^\gamma}, & \text{if } 0 < \gamma \le 1, \end{cases} \quad \text{a.s.,}$$

whenever $\lambda \to 0$ and $r \gg \log(1/\lambda)\log\log\log(1/\lambda)$.

Let us admit (3.4) for the time being, and prove the lemma. By the occupation time formula and a change of variables, we have

$$\Phi(\tau(r)) = \int_{-\infty}^\infty \exp\left(-\frac{4}{S^{-1}(x)} - 8\lambda S^{-1}(x)\right) L(\tau(r), x)\, dx,$$

which means that $\Phi(\tau(r)) = D_0(r)$ in the notation of (3.3). By (3.4), whenever $\lambda \to 0$ and $r \gg \log(1/\lambda)\log\log\log(1/\lambda)$,

$$\Phi(\tau(r)) = (1+o(1))r\log(1/\lambda) \quad \text{a.s.,}$$

which means that whenever $\lambda \to 0$ and $v \gg \log^2(1/\lambda)\log\log\log(1/\lambda)$,

$$(3.5) \quad \Phi^{-1}(v) = \tau\left((1+o(1))\frac{v}{\log(1/\lambda)}\right) \quad \text{a.s.}$$

On the other hand, by (3.2) and the occupation time formula,

$$\begin{aligned}
\int_0^v Z(s)\, ds &= \int_0^v S^{-1}(B(\Phi^{-1}(s)))\, ds \\
&= \int_0^{\Phi^{-1}(v)} \frac{S^{-1}(B(r))}{h^2(B(r))}\, dr \\
&= \int_{-\infty}^\infty \frac{S^{-1}(x)}{h^2(x)} L(\Phi^{-1}(v), x)\, dx \\
&\stackrel{\text{def}}{=} \Psi(\Phi^{-1}(v)),
\end{aligned}$$

with obvious definition of $\Psi(\cdot)$. Note that $\Psi(\tau(r)) = D_1(r)$ in the notation of (3.3). According to (3.4),

$$\Psi(\tau(r)) = (1+o(1))\frac{r}{4\lambda}, \quad r \gg \log(1/\lambda)\log\log\log(1/\lambda) \quad \text{a.s.}$$

which, with the aid of (3.5), would yield Lemma 3.2.

It remains to show (3.4). We shall make use of the following simple consequence of law of large numbers: if $f:\mathbb{R}\to\mathbb{R}$ is such that $\int_\mathbb{R}|f(x)|\,dx < \infty$, then

$$(3.6) \quad \lim_{r\to\infty} \frac{1}{r}\int_\mathbb{R} f(x) L(\tau(r), x)\, dx = \int_\mathbb{R} f(x)\, dx \quad \text{a.s.}$$

We write
$$D_\gamma(r) = \left(\int_{-\infty}^0 + \int_0^{S(a)} + \int_{S(a)}^\infty\right)(S^{-1}(x))^\gamma$$
$$\times \exp\left(-\frac{4}{S^{-1}(x)} - 8\lambda S^{-1}(x)\right)L(\tau(r),x)\,dx$$
$$\stackrel{\text{def}}{=} \Delta_1 + \Delta_2 + \Delta_3,$$

where $a = a(r,\lambda) > 1/\lambda$ is chosen such that

(3.7)
$$\frac{e^{4\lambda a}}{4\lambda a} = \log\left(\frac{1}{\lambda}\right)\log\log\left(\frac{1}{\lambda}\right) \quad \text{if } r \geq \log^2\left(\frac{1}{\lambda}\right),$$
$$a = \frac{1}{4\lambda}\log\log\left(\frac{1}{\lambda}\right) \quad \text{if } r < \log^2\left(\frac{1}{\lambda}\right).$$

For $x < 0$, let $S^{-1}(x) = y \in (0,1)$, so that for $\lambda < \frac{1}{4}$, $-x = \int_y^1 z^{-1}e^{4\lambda z + 2/z}\,dz \leq y^{-1}e^{1+2/y}$. This implies the existence of a large constant $c > 0$ such that for all $x < -c$ and $\lambda < \frac{1}{4}$, we have $\frac{2}{S^{-1}(x)} \geq \log|x| - 2\log\log|x|$. Hence, by means of (3.6),

(3.8)
$$\Delta_1 \leq \int_{-\infty}^0 e^{-8\lambda S^{-1}(x)}L(\tau(r),x)\,dx$$
$$\leq \int_{-\infty}^0 \left(\mathbf{1}_{(x \geq -c)} + \mathbf{1}_{(x < -c)}\frac{\log^4|x|}{|x|^2}\right)L(\tau(r),x)\,dx = O(r) \quad \text{a.s.}$$

To treat $\Delta_3$, we observe that for $y \geq a$ (thus, $y\lambda \to \infty$),

(3.9)
$$S(y) = \int_1^y e^{(2/z)+4\lambda z}\,\frac{dz}{z}$$
$$= (1+o(1))\int_1^y e^{4\lambda z}\,\frac{dz}{z}$$
$$= (1+o(1))\left(\int_{4\lambda}^1 e^x\,\frac{dx}{x} + \int_1^{4y\lambda} e^x\,\frac{dx}{x}\right)$$
$$= (1+o(1))\left(\log\left(\frac{1}{\lambda}\right) + O(1) + \frac{1+o(1)}{4\lambda y}e^{4\lambda y}\right).$$

Let us distinguish two cases: First, if $r \geq \log^2(1/\lambda)$, then $\frac{1}{4a\lambda}e^{4a\lambda} \gg \log(1/\lambda)$; hence, for $y \geq a$,

(3.10)
$$S(y) = (1+o(1))\frac{e^{4\lambda y}}{4\lambda y}.$$

Thus, for $x \geq S(a)$, we have
$$e^{4\lambda S^{-1}(x)} = (1+o(1))4\lambda S^{-1}(x)x > x \quad \text{and} \quad S^{-1}(x) \sim \frac{\log x}{4\lambda}.$$

It follows from (3.6) that

$$\Delta_3 \leq \int_{S(a)}^{\infty} (S^{-1}(x))^{\gamma} e^{-8\lambda S^{-1}(x)} L(\tau(r), x) \, dx$$

$$\leq 2 \int_{S(a)}^{\infty} \left(\frac{\log x}{4\lambda}\right)^{\gamma} x^{-2} L(\tau(r), x) \, dx$$

$$\leq 2(S(a))^{-1/2} \int_{1}^{\infty} \left(\frac{\log x}{4\lambda}\right)^{\gamma} x^{-3/2} L(\tau(r), x) \, dx$$

$$= o\left(\frac{r}{\lambda^{\gamma}}\right).$$

For the other case $\log(1/\lambda) \log \log \log(1/\lambda) \ll r < \log^2(1/\lambda)$, we have

$$\Delta_3 \leq L^*(\tau(r)) \int_{S(a)}^{\infty} (S^{-1}(x))^{\gamma} \exp\left\{-\frac{4}{S^{-1}(x)} - 8S^{-1}(x)\right\} dx$$

$$= L^*(\tau(r)) \int_{a}^{\infty} e^{-(2/y)-4\lambda y} \frac{dy}{y^{1-\gamma}}$$

$$\leq L^*(\tau(r)) \int_{a}^{\infty} e^{-4\lambda y} \, dy = \frac{L^*(\tau(r))}{4\lambda a} < L^*(\tau(r)),$$

where

$$L^*(\tau(r)) \stackrel{\text{def}}{=} \sup_{x \in \mathbb{R}} L(\tau(r), x) \leq r(\log r)^{1+\varepsilon} \qquad \text{a.s.,}$$

for any $\varepsilon > 0$, by the law of iterated logarithm for $L^*(\cdot)$ ([14]). Since $r < \log^2(1/\lambda)$, this implies that $\Delta_3 \ll r(\log \log(1/\lambda))^2$, a.s. Therefore, in both situations, we have

(3.11) $$\Delta_3 = o\left(\frac{r}{\lambda^{\gamma}}\right) + o\left(r\left(\log \log\left(\frac{1}{\lambda}\right)\right)^2\right) \qquad \text{a.s.}$$

To deal with $\Delta_2$, we claim that in both cases [i.e., $r \geq \log^2(1/\lambda)$ and $r < \log^2(1/\lambda)$],

$$S(a) = o\left(\frac{r}{\log \log r}\right).$$

In fact, if $r \geq \log^2(1/\lambda)$, then by (3.10), $S(a) \leq \frac{2}{4a\lambda} e^{4a\lambda} = 2\log(1/\lambda) \log \log(1/\lambda) = o(r/\log \log r)$; otherwise, by definition (3.7), $a = \frac{\log \log(1/\lambda)}{4\lambda}$ and by means of (3.9), we have $S(a) \leq 2\log(1/\lambda) = o(r/\log \log r)$ [recalling that $r \gg \log(1/\lambda) \log \log \log(1/\lambda)$]. Hence, we can apply Lemma 2.4 to see that

$$\Delta_2 = r(1+o(1)) \int_0^{S(a)} (S^{-1}(x))^{\gamma} \exp\left(-\frac{4}{S^{-1}(x)} - 8\lambda S^{-1}(x)\right) dx \qquad \text{a.s.}$$

By a change of variables $S^{-1}(x) = y$, the integral on the right-hand side is, when $\lambda \to 0$,

$$\begin{aligned}
&= \int_1^a e^{-(2/y)-4\lambda y} \frac{dy}{y^{1-\gamma}} = (1+o(1)) \int_1^a e^{-4\lambda y} \frac{dy}{y^{1-\gamma}} \\
&= \frac{1+o(1)}{(4\lambda)^\gamma} \int_{4\lambda}^{4\lambda a} e^{-u} \frac{du}{u^{1-\gamma}} \\
&= \begin{cases} (1+o(1))\log(1/\lambda), & \text{if } \gamma = 0, \\ (1+o(1))\frac{\Gamma(\gamma)}{(4\lambda)^\gamma}, & \text{if } \gamma > 0. \end{cases}
\end{aligned}$$

Accordingly,

$$(3.12) \quad \Delta_2 = \begin{cases} (1+o(1))r\log(1/\lambda), & \text{if } \gamma = 0, \\ (1+o(1))\Gamma(\gamma)\frac{r}{(4\lambda)^\gamma}, & \text{if } \gamma > 0, \end{cases} \quad \text{a.s.}$$

Since $D_\gamma(r) = \Delta_1 + \Delta_2 + \Delta_3$, assembling (3.8), (3.11) and (3.12) readily yields (3.4). The lemma is proved. □

We now have all the ingredients for the proof of Theorem 1.1.

PROOF OF THEOREM 1.1.  Let $t, v \to \infty$ satisfy the conditions in Theorem 1.1.

First, we prove (1.3) for $\sup_{0 \le s \le t} X(s)$ in place of $X(t)$: $P$-almost surely,

$$\mathbf{P}_\omega\left(\sup_{0 \le s \le t} X(s) > v\right) = \mathbf{P}_\omega(H(v) < t) = \exp\left(-(1+o(1))\frac{v}{2\log(t/v)}\right).$$
(3.13)

To this end, applying Chebyshev's inequality to Fact 3.1 and Lemma 3.2, we have that for almost surely all potentials $W$ and for $\lambda \to 0$ satisfying $v \gg \log^2(1/\lambda) \log\log\log(1/\lambda)$,

$$(3.14) \quad \mathbf{P}_\omega\{H(v) < t\} \le e^{\lambda t} \mathbf{E}_\omega[e^{-\lambda H(v)}] = \exp\left(\lambda t - (1+o(1))\frac{v}{2\log(1/\lambda)}\right).$$

By choosing $\lambda = \frac{v}{t \log^2(t/v)}$ [this is possible since $t \gg v \gg \log^2(t/v) \log\log\log(t/v)$], we have

$$\mathbf{P}_\omega\{H(v) < t\} \le \exp\left(-(1+o(1))\frac{v}{2\log(t/v)}\right) \quad \text{a.s.}$$

This implies the upper bound in (3.13).

To get the lower bound, we keep the choice of $\lambda$ and use the simple relation

$$\mathbf{P}_\omega\{H(v) < t\log^2 t\} \geq \mathbf{E}_\omega[e^{-\lambda H(v)}] - e^{-\lambda t \log^2 t}$$
(3.15)
$$\geq \exp\left(-(1+o(1))\frac{v}{2\log(t/v)}\right),$$

by means of Lemma 3.2. Write $\tilde{t} = t \log^2 t$. Since $\log \log t = o(\log(t/v))$, we have $\log(\tilde{t}/v) \sim \log(t/v)$. Thus, (3.15) yields the lower bound in (3.13).

In light of the trivial inequality $\mathbf{P}_\omega\{X(t) > v\} \leq \mathbf{P}_\omega(H(v) < t)$, it remains to show the lower bound in (1.3). For any $c_1, c_2 > 0$, we define

$$G(c_1, c_2) \stackrel{\text{def}}{=} \{(v,t) : t, v \geq 3^{81}, v \geq c_1 \log^2 t \log \log \log t, \log(t/v) \geq c_2 \log \log t\}.$$

The proof of Theorem 1.1 will be complete if we can show that for any small $\varepsilon > 0$ and almost surely all $\omega$, there exist $c_1(\varepsilon, \omega), c_2(\varepsilon, \omega) > 0$ and $t_0 = t_0(\varepsilon, \omega) > 0$ such that for all $t, v \geq t_0$ and $(v,t) \in G(c_1(\varepsilon, \omega), c_2(\varepsilon, \omega))$, we have

(3.16) $$\log \mathbf{P}_\omega\{X(t) > (1-\varepsilon)v\} \geq -\frac{1+2\varepsilon}{2}\frac{v}{\log(t/v)}.$$

Write $\mathbf{P}_{x,\omega}$ for the law of the diffusion $X$ defined in (1.1) starting from $X(0) = x$. Let $j$, $k \geq j_0$ be sufficiently large and define $v_j = e^{j/\log j}$ and $t_k = \exp(e^{k/\log k})$. We claim that for $c_3 = c_3(\varepsilon) \stackrel{\text{def}}{=} \frac{256}{(1-2\varepsilon)\pi^2 \varepsilon}$,

(3.17) $$\sum_{j,k \geq j_0, (v_j, t_k) \in G(c_3, 1)} P\{\mathbf{P}_{v_j,\omega}\{H((1-\varepsilon)v_j) < t_k\} > \tfrac{1}{2}\} < \infty.$$

Indeed, by Brownian symmetry, $\mathbf{P}_{v_j,\omega}\{H((1-\varepsilon)v_j) < t_k\}$ is distributed as $\mathbf{P}_{0,\omega}\{H(\varepsilon v_j) < t_k\}$, so that by means of (5.8) in Section 5, for all large $j$ and $k$,

(3.18)
$$\begin{aligned}
&E(\mathbf{P}_{v_j,\omega}\{H((1-\varepsilon)v_j) < t_k\}) \\
&= \mathbb{P}\{H(\varepsilon v_j) < t_k\} \\
&\leq c\varepsilon^2 v_j^2 e^{-\varepsilon v_j} + 9\exp\left(-(1-2\varepsilon)\frac{\pi^2}{8}\frac{\varepsilon v_j}{\log^2(\varepsilon v_j t_k)}\right) \\
&\leq (c+9)\exp\left(-\frac{(1-2\varepsilon)\pi^2 \varepsilon}{32}\frac{v_j}{\log^2 t_k}\right).
\end{aligned}$$

We have used in the last inequality the fact that $v_j \leq t_k$.

By Chebyshev's inequality,

$$\sum_{j,k \geq j_0, (v_j, t_k) \in G(c_3, 1)} P\left\{\mathbf{P}_{v_j,\omega}\{H((1-\varepsilon)v_j) < t_k\} > \frac{1}{2}\right\}$$

$$\leq 2 \sum_{j,k \geq j_0, (v_j, t_k) \in G(c_3, 1)} E(\mathbf{P}_{v_j,\omega}\{H((1-\varepsilon)v_j) < t_k\})$$

$$\leq 2(c+9) \sum_{j,k \geq j_0, (v_j,t_k) \in G(c_3,1)} \exp\left(-\frac{(1-2\varepsilon)\pi^2 \varepsilon}{32} \frac{v_j}{\log^2 t_k}\right).$$

Several elementary computations show that the above (double) sum is finite: in fact, we can decompose this sum into $\sum_{j \geq k^4}$ and $\sum_{j \leq k^4}$. Note that for $j \geq k^4$, $\frac{v_j}{\log^2 t_k} = \exp\{\frac{j}{\log j} - 2\frac{k}{\log k}\} \geq \exp\{\frac{j}{2\log j}\}$. Hence, $\sum_{j \geq k^4} \leq \sum_j j^{1/4} \times \exp(-\frac{(1-2\varepsilon)\pi^2 \varepsilon}{32} \exp\{\frac{j}{2\log j}\}) < \infty$. For the case $j \leq k^4$, we use the definition of $G(c_3, 1)$, which says that $v_j/\log^2 t_k \geq \frac{256}{(1-2\varepsilon)\pi^2 \varepsilon} \log\log\log t_k$, $\sum_{j \leq k^4} \leq \sum_k k^4 e^{-8\log\log\log t_k} \leq \sum_k k^{-4} < \infty$. This yields (3.17).

We now proceed to the proof of (3.16). According to the lower bound in (3.13), for $P$-almost all $\omega$, there exist $c_1(\varepsilon, \omega) \geq 2c_3$ [recalling that $c_3 = \frac{256}{(1-2\varepsilon)\pi^2 \varepsilon}$] and $c_2(\varepsilon, \omega) \geq 2$ such that for all $(v,t) \in G(\frac{c_1(\varepsilon,\omega)}{2}, \frac{c_2(\varepsilon,\omega)}{2})$,

$$(3.19) \qquad \mathbf{P}_{0,\omega}\{H(v) < t\} \geq \exp\left(-\frac{1+\varepsilon}{2\log(t/v)}\right).$$

For any $\omega$, the strong Markov property at time $H(v_j)$ implies that for any $v_{j-1} \leq v \leq v_j$ and $t_{k-1} \leq t \leq t_k$,

$$\mathbf{P}_{0,\omega}\{X(t) > (1-\varepsilon)v\} \geq \mathbf{P}_{0,\omega}\{H(v_j) < t_{k-1}\}\mathbf{P}_{v_j,\omega}\{H((1-\varepsilon)v_j) \geq t_k\}.$$
(3.20)

Consider large $t, v$ such that $(v,t) \in G(c_1(\varepsilon, \omega), c_2(\varepsilon, \omega))$, say $v \in [v_{j-1}, v_j]$ and $t \in [t_{k-1}, t_k]$. Then both $(v_j, t_{k-1})$ and $(v_j, t_k)$ are elements of $G(\frac{c_1(\varepsilon,\omega)}{2}, \frac{c_2(\varepsilon,\omega)}{2})$, and a fortiori $(v_j, t_k) \in G(c_3, 1)$. According to (3.19),

$$\mathbf{P}_{0,\omega}\{H(v_j) < t_{k-1}\} \geq \exp\left(-\frac{1+\varepsilon}{2\log(t_{k-1}/v_j)}\right),$$

whereas by (3.17) and the Borel–Cantelli lemma,

$$\mathbf{P}_{v_j,\omega}\{H((1-\varepsilon)v_j) \geq t_k\} \geq \tfrac{1}{2}.$$

Plugging these estimates into (3.20) yields

$$\mathbf{P}_{0,\omega}\{X(t) > (1-\varepsilon)v\} \geq \frac{1}{2}\exp\left(-\frac{(1+\varepsilon)v_j}{2\log(t_{k-1}/v_j)}\right) \geq \exp\left(-\frac{(1+2\varepsilon)v}{2\log(t/v)}\right),$$

since $v_j/v_{j-1} \to 1$ and $\log t_k / \log t_{k-1} \to 1$. This will yield the lower bound (3.16). □

**4. The annealed case: a joint arcsine law.** Let $\kappa \in \mathbb{R}$ and let

$$W_\kappa(x) = W(x) - \frac{\kappa}{2}x, \qquad x \in \mathbb{R},$$

where $W$ is, as before, a Brownian motion defined on $\mathbb{R}$ with $W(0) = 0$. In this section we shall study the diffusion $X$ with potential $W_\kappa$ [i.e., replacing $W$ by $W_\kappa$ in (1.1)]. Plainly, when $\kappa = 0$, we recover the case of Brownian potential and $X$ is recurrent, whereas $X(t) \to +\infty$, $\mathbb{P}$-a.s. if $\kappa > 0$.

We recall the time change representation of $X$ (cf. [2] for $\kappa = 0$ and [10] for $\kappa > 0$):

$$X(t) = A_\kappa^{-1}(B(T_\kappa^{-1}(t))), \qquad t \geq 0,$$

where $B$ is a one-dimensional Brownian motion starting from 0, independent of $W$, and

$$A_\kappa(x) = \int_0^x e^{W_\kappa(y)} \, dy, \qquad x \in \mathbb{R},$$

$$T_\kappa(t) = \int_0^t e^{-2W_\kappa(A_\kappa^{-1}(B(s)))} \, ds, \qquad t \geq 0.$$

(Recall that $A_\kappa^{-1}$ and $T_\kappa^{-1}$ denote the inverses of $A_\kappa$ and $T_\kappa$, resp.) Here, we stress the fact that the process $B$, a Brownian motion independent of $W$, is not the same $B$ as in Section 3. There is no risk of confusion since we always separate the quenched and the annealed cases. Recall that $(L(t,x), t \geq 0, x \in \mathbb{R})$ denote the local times of $B$ and $\sigma(\cdot)$ is the process of first hitting times of $B$:

$$\sigma(x) \stackrel{\text{def}}{=} \inf\{t > 0 : B(t) > x\}, \qquad x \geq 0.$$

Therefore, $H(v)$, the first hitting time of $X$ at $v > 0$, can be represented as follows:

$$H(v) = \inf\{t \geq 0 : X(t) = v\}$$
$$(4.1) \qquad = T_\kappa(\sigma(A_\kappa(v)))$$
$$= \int_0^{\sigma(A_\kappa(v))} e^{-2W_\kappa(A_\kappa^{-1}(B(s)))} \, ds$$
$$= \int_{-\infty}^{A_\kappa(v)} e^{-2W_\kappa(A_\kappa^{-1}(x))} L(\sigma(A_\kappa(v)), x) \, dx$$
$$(4.2) \qquad = \Theta_1(v) + \Theta_2(v),$$

with

$$\Theta_1(v) = \int_0^{A_\kappa(v)} e^{-2W_\kappa(A_\kappa^{-1}(x))} L(\sigma(A_\kappa(v)), x) \, dx = \int_0^{H(v)} \mathbf{1}_{\{X(s) \geq 0\}} \, ds,$$

$$\Theta_2(v) = \int_{-\infty}^0 e^{-2W_\kappa(A_\kappa^{-1}(x))} L(\sigma(A_\kappa(v)), x) \, dx = \int_0^{H(v)} \mathbf{1}_{\{X(s) < 0\}} \, ds.$$

The main result in this section describes the (annealed) distribution of $(\Theta_1(v), \Theta_2(v))$.

THEOREM 4.1. *Let $\kappa \geq 0$ and let $v > 0$. Under $\mathbb{P}$, we have*

$$(\Theta_1(v), \Theta_2(v)) \stackrel{\text{law}}{=} \left(4\int_0^v (e^{\Xi_\kappa(s)} - 1)\, ds, 16\Upsilon_{2-2\kappa}(e^{\Xi_\kappa(v)/2} \leadsto 1)\right),$$

*where $\Upsilon_{2-2\kappa}(x \leadsto y)$ denotes the first hitting time of $y$ by a $\mathrm{BES}(2-2\kappa)$ starting from $x$, independent of the diffusion $\Xi_\kappa$ which is the unique non-negative solution of*

$$(4.3)\ \Xi_\kappa(t) = \int_0^t \sqrt{1 - e^{-\Xi_\kappa(s)}}\, d\beta(s) + \int_0^t \left(-\frac{\kappa}{2} + \frac{1+\kappa}{2} e^{-\Xi_\kappa(s)}\right) ds, \qquad t \geq 0,$$

*$\beta$ being a standard Brownian motion.*

The proof of Theorem 4.1 involves some deep results. Let us first recall Lamperti's representation theorem for exponential functionals [15].

FACT 4.2 (Lamperti's representation). Let $\kappa \in \mathbb{R}$. There exists a $\mathrm{BES}(2+2\kappa)$, denoted by $\widetilde{R}$ and starting from $\widetilde{R}(0) = 2$, such that

$$(4.4)\quad \exp\left(W(x) + \frac{\kappa}{2}x\right) = \frac{1}{4}\widetilde{R}^2\left(\int_0^x \exp\left\{W(y) + \frac{\kappa y}{2}\right\} dy\right), \qquad x \geq 0.$$

Let $d_1, d_2 \geq 0$, $a \in [0,1]$, and consider the equation

$$(4.5)\quad\begin{aligned} dY(t) &= 2\sqrt{Y(t)(1 - Y(t))}\, d\beta(t) + (d_1 - (d_1 + d_2)Y(t))\, dt,\\ Y(0) &= a,\end{aligned}$$

where $\beta$ is a standard one-dimensional Brownian motion. The solution $Y$ of the above equation is called a Jacobi process of dimension $(d_1, d_2)$, starting from $a$ (see [11]). We mention that almost surely, $0 \leq Y(t) \leq 1$ for all $t \geq 0$.

The following result gives the skew-product representation of two independent Bessel processes in terms of the Jacobi process.

FACT 4.3 ([25]). Let $R_1$ and $R_2$ be two independent Bessel processes of dimensions $d_1$ and $d_2$, respectively. We assume $d_1 + d_2 \geq 2$, $R_1(0) = r_1 \geq 0$ and $R_2(0) = r_2 > 0$. Then there exists a Jacobi process $Y$ of dimension $(d_1, d_2)$ starting from $\frac{r_1^2}{r_1^2 + r_2^2}$, independent of the process $(R_1^2(t) + R_2^2(t), t \geq 0)$, such that

$$\frac{R_1^2(t)}{R_1^2(t) + R_2^2(t)} = Y\left(\int_0^t \frac{ds}{R_1^2(s) + R_2^2(s)}\right), \qquad t \geq 0.$$

We are now ready to prove Theorem 4.1.

PROOF OF THEOREM 4.1. Using Fact 2.1 and the independence of $B$ and $W$, the process $\{\frac{1}{A_\kappa(v)} L(\sigma(A_\kappa(v)), (1-x)A_\kappa(v)), x \geq 0\}$ is a strong

Markov process starting from 0, independent of $W$; it is a BESQ(2) for $x \in [0, 1]$, and is a BESQ(0) for $x \geq 1$. By scaling,

$$(\Theta_1(v), \Theta_2(v)) \stackrel{\text{law}}{=} \left( \int_0^v e^{-W_\kappa(x)} R^2(A_\kappa(v) - A_\kappa(x)) \, dx, \int_{-\infty}^0 e^{-2W_\kappa(A_\kappa^{-1}(x))} U(|x|) \, dx \right),$$

where $R$ is a BES(2) starting from 0, independent of $W$, and conditionally on $(R, W)$, $U$ is a BESQ(0) starting from $R^2(A_\kappa(v))$.

By time reversal, $(\widehat{W}_\kappa(y) \stackrel{\text{def}}{=} W_\kappa(v - y) - W_\kappa(v), 0 \leq y \leq v)$ has the same law as $(W_{-\kappa}(y) = W(y) + \frac{\kappa}{2} y, 0 \leq y \leq v)$. Observe that

$$\int_0^v e^{-W_\kappa(x)} R^2(A_\kappa(v) - A_\kappa(x)) \, dx$$
$$= \int_0^v e^{-\widehat{W}_\kappa(y)} e^{\widehat{W}_\kappa(v)} R^2 \left( e^{-\widehat{W}_\kappa(v)} \int_0^y e^{\widehat{W}_\kappa(z)} \, dz \right) dy,$$
$$R^2(A_\kappa(v)) = R^2 \left( e^{-\widehat{W}_\kappa(v)} \int_0^v e^{\widehat{W}_\kappa(z)} \, dz \right).$$

By scaling and independence of $R$ and $W$, the process $x \mapsto e^{\widehat{W}_\kappa(v)} R^2(x e^{-\widehat{W}_\kappa(v)})$ has the same law as $R$ and is independent of $W$. It follows that

$$(\Theta_1(v), \Theta_2(v))$$
(4.6)
$$\stackrel{\text{law}}{=} \left( \int_0^v e^{-W_{-\kappa}(x)} R^2(A_{-\kappa}(x)) \, dx, \int_{-\infty}^0 e^{-2W_\kappa(A_\kappa^{-1}(x))} U(|x|) \, dx \right),$$

where conditionally on $(R, W)$, $U$ has the same law as a BESQ(0) starting from $U(0) = e^{-W_{-\kappa}(v)} R^2(A_{-\kappa}(v))$.

Let us first treat the part $(W_{-\kappa}(x), x \geq 0)$. By means of Fact 4.2, there exists a Bessel process $\widetilde{R}$ of dimension $2 + 2\kappa \geq 2$, starting from 2, such that

(4.7)
$$\int_0^v e^{-W_{-\kappa}(x)} R^2(A_{-\kappa}(x)) \, dx = \int_0^{A_{-\kappa}(v)} e^{-2W_{-\kappa}(A_{-\kappa}^{-1}(y))} R^2(y) \, dy$$
$$= 16 \int_0^{A_{-\kappa}(v)} \frac{R^2(y)}{\widetilde{R}^4(y)} \, dy,$$

where we stress the independence of the two Bessel processes $R$ and $\widetilde{R}$. Observe that

$$A_{-\kappa}^{-1}(x) = 4 \int_0^x \frac{du}{\widetilde{R}^2(u)}, \qquad x \geq 0.$$

We apply Fact 4.3 to $R$ and $\widetilde{R}$, to see that there exists a Jacobi process $Y$ of dimension $(2, 2 + 2\kappa)$ starting from 0 such that

$$\frac{R^2(x)}{R^2(x) + \widetilde{R}^2(x)} = Y \left( \int_0^x \frac{ds}{R^2(s) + \widetilde{R}^2(s)} \right) \stackrel{\text{def}}{=} Y(\Lambda(x)), \qquad x \geq 0,$$

where $\Lambda(x) \stackrel{\text{def}}{=} \int_0^x \frac{ds}{R^2(s)+\widetilde{R}^2(s)}$, and $\Lambda(\cdot)$ is independent of $Y$. Note that $Y(0) = 0$ and $0 < Y(t) < 1$ for all $t > 0$.

This representation, together with (4.7), implies that

$$\int_0^v e^{-W_{-\kappa}(x)} R^2(A_{-\kappa}(x))\, dx = 16 \int_0^{\Lambda(A_{-\kappa}(v))} \frac{Y(u)}{(1-Y(u))^2}\, du$$
(4.8)
$$= 16 \int_0^{\rho^{-1}(v)} \frac{Y(u)}{(1-Y(u))^2}\, du,$$

where

$$\rho(x) \stackrel{\text{def}}{=} A_{-\kappa}^{-1}(\Lambda^{-1}(x)) = 4 \int_0^{\Lambda^{-1}(x)} \frac{du}{\widetilde{R}^2(u)} = 4 \int_0^x \frac{dy}{1-Y(y)},$$

by a change of variables $y = \Lambda(u)$. Going back to (4.6),

$$U(0) = e^{-W_{-\kappa}(v)} R^2(A_{-\kappa}(v))$$
$$= \frac{4R^2(A_{-\kappa}(v))}{\widetilde{R}^2(A_{-\kappa}(v))} = \frac{4Y(\Lambda(A_{-\kappa}(v)))}{1-Y(\Lambda(A_{-\kappa}(v)))} = \frac{4Y(\rho^{-1}(v))}{1-Y(\rho^{-1}(v))}.$$

Assume for the moment that for any fixed $r > 0$, if $U_r$ denotes a BESQ(0) starting from $U_r(0) = r$, independent of $W$, then

$$(4.9) \quad \int_{-\infty}^0 e^{-2W_\kappa(A_\kappa^{-1}(x))} U_r(|x|)\, dx \stackrel{\text{law}}{=} 16 \Upsilon_{2-2\kappa}\left(\frac{\sqrt{4+r}}{2} \rightsquigarrow 1\right), \qquad \kappa \geq 0.$$

By admitting (4.9), it follows from (4.6), (4.7) and (4.8) that under the total probability $\mathbb{P}$,

(4.10)
$$(\Theta_1(v), \Theta_2(v))$$
$$\stackrel{\text{law}}{=} \left(16 \int_0^{\rho^{-1}(v)} \frac{Y(u)}{(1-Y(u))^2}\, du, 16 \Upsilon_{2-2\kappa}\left(\frac{\sqrt{4+U(0)}}{2} \rightsquigarrow 1\right)\right),$$

where $U(0) = \frac{4Y(\rho^{-1}(v))}{1-Y(\rho^{-1}(v))}$, and given $U(0) = r$, $\Upsilon_{2-2\kappa}(\frac{\sqrt{4+r}}{2} \rightsquigarrow 1)$ is the first hitting time of 1 by a BES$(2-2\kappa)$ starting from $\frac{\sqrt{4+r}}{2}$, independent of the process $Y$.

Since $d\rho^{-1}(x) = \frac{1-Y(\rho^{-1}(x))}{4}\, dx$ and

$$\int_0^{\rho^{-1}(v)} \frac{Y(u)}{(1-Y(u))^2}\, du = \frac{1}{4} \int_0^v \left(\frac{1}{1-Y(\rho^{-1}(x))} - 1\right) dx, \qquad v > 0,$$

it follows from (4.5) that $\Xi_\kappa(t) \stackrel{\text{def}}{=} -\log\{1-Y(\rho^{-1}(t))\}$ satisfies the stochastic integral equation (4.3). Theorem 4.1 will then follow from the identity in law (4.10).

It remains to show (4.9). Note that $(W_\kappa(-x), x \geq 0)$ is distributed as $(W_{-\kappa}(x), x \geq 0)$. Thus,

$$
\text{(4.11)} \quad \int_{-\infty}^{0} e^{-2W_\kappa(A_\kappa^{-1}(x))} U_r(|x|) \, dx \stackrel{\text{law}}{=} \int_{0}^{\infty} e^{-2W_{-\kappa}(A_{-\kappa}^{-1}(x))} U_r(x) \, dx
$$
$$
= 16 \int_{0}^{\infty} \frac{U_r(x)}{\widetilde{R}^4(x)} \, dx,
$$

by using again the BES$(2+2\kappa)$ process $\widetilde{R}$ defined in (4.7).

Applying Fact 4.3 to the two independent squared Bessel processes $U_r$ and $\widetilde{R}^2$, we get a Jacobi process $\widehat{Y}$ of dimension $(0, 2+2\kappa)$ starting from $\frac{r}{r+4}$, such that

$$
\frac{U_r(t)}{U_r(t) + \widetilde{R}^2(t)} = \widehat{Y}\left(\int_0^t \frac{ds}{U_r(s) + \widetilde{R}^2(s)}\right) \stackrel{\text{def}}{=} \widehat{Y}(\widehat{\Lambda}(t)), \quad t \geq 0,
$$

with $\widehat{\Lambda}(t) \stackrel{\text{def}}{=} \int_0^t \frac{ds}{U_r(s) + \widetilde{R}^2(s)}$, $t \geq 0$, independent of $\widehat{Y}$. Observe that $\widehat{Y}$ is absorbed at 0.

By a change of variables $t = \widehat{\Lambda}(x)$,

$$
\text{(4.12)} \quad 16 \int_0^{\infty} \frac{U_r(x)}{\widetilde{R}^4(x)} \, dx = 16 \int_0^{T_{\widehat{Y}}(0)} \frac{\widehat{Y}(t)}{(1-\widehat{Y}(t))^2} \, dt,
$$

where $T_{\widehat{Y}}(0) \stackrel{\text{def}}{=} \inf\{t : \widehat{Y}(t) = 0\}$. By computing the scale function and using the Dubins–Schwarz theorem ([19], Theorem V.1.6) for continuous local martingales, there exists some one-dimensional Brownian motion $\beta$ starting from 0 such that

$$
s(\widehat{Y}(t)) = \beta(\phi(t)), \quad t \geq 0,
$$

with

$$
y_0 \stackrel{\text{def}}{=} \frac{r}{r+4},
$$

$$
s(y) \stackrel{\text{def}}{=} \begin{cases} \log \dfrac{1-y_0}{1-y}, & \text{if } \kappa = 0, \\ \dfrac{1}{\kappa}\{(1-y)^{-\kappa} - (1-y_0)^{-\kappa}\}, & \text{if } \kappa > 0, \end{cases} \quad 0 \leq y < 1,
$$

$$
\phi(t) \stackrel{\text{def}}{=} 4 \int_0^t \frac{\widehat{Y}(s)}{(1-\widehat{Y}(s))^{1+2\kappa}} \, ds, \quad t \geq 0.
$$

Note that $\phi(T_{\widehat{Y}}(0)) = \inf\{t > 0 : \beta(t) = s(0)\} \stackrel{\text{def}}{=} T_\beta(s(0))$.

When $\kappa = 0$, we have

$$16 \int_0^{T_{\widehat{Y}}(0)} \frac{\widehat{Y}(t)}{(1-\widehat{Y}(t))^2} dt$$

$$= \frac{4}{1-y_0} \int_0^{T_\beta(s(0))} e^{\beta(u)} du$$

$$\stackrel{\text{law}}{=} \frac{4}{1-y_0} \inf\{s > 0 : \widetilde{R}(s) = 2e^{s(0)/2} = 2\sqrt{1-y_0}\},$$

where the last equality in law follows from (4.4) by replacing $W$ by $\beta$ [recalling $\widetilde{R}(0) = 2$]. This, together with the scaling property, yields (4.9) in the case $\kappa = 0$.

If $\kappa > 0$, we observe that

$$16 \int_0^{T_{\widehat{Y}}(0)} \frac{\widehat{Y}(t)}{(1-\widehat{Y}(t))^2} dt = 4 \int_0^{T_\beta(s(0))} (1 - s^{-1}(\beta(u)))^{2\kappa - 1} du$$

$$= 4\kappa^{(1/\kappa)-2} \int_0^{T_\beta(s(0))} \left(\frac{1}{\kappa} - s(0) + \beta(u)\right)^{(1/\kappa)-2} du.$$

By symmetry (i.e., replacing $\beta$ by $-\beta$) and Lemma 2.3, the expression on the right-hand side is equal in law to

$$4\kappa^{(1/\kappa)-2} \int_0^{T_\beta(|s(0)|)} \left(\frac{1}{\kappa} + |s(0)| - \beta(u)\right)^{(1/\kappa)-2} du \stackrel{\text{law}}{=} 16\Upsilon_{2-2\kappa}\left(\frac{\sqrt{4+r}}{2} \rightsquigarrow 1\right),$$

completing the proof of (4.9). Theorem 4.1 is proved. $\square$

**5. Proof of Theorem 1.2.** This section is devoted to the proof of Theorem 1.2. We prove the upper and lower bounds with different approaches.

5.1. *Theorem* 1.2: *the upper bound.* The proof of the upper bound is based on an analysis of the diffusion process $\Xi_0(\cdot)$ introduced in (4.3) (with $\kappa = 0$).

For notational convenience, we write $\Xi \stackrel{\text{def}}{=} \Xi_0$. Let us start with a couple of lemmas.

LEMMA 5.1. *There exists a numerical constant $c > 0$ such that for all $t > 100$ and $0 < a < x < \sqrt{t}$, we have*

$$\mathbb{P}\left(\sup_{0 \leq t_1 \leq t_2 < t, t_2 - t_1 < a} |\Xi(t_2) - \Xi(t_1)| > x\right) \leq c\frac{t}{a} \exp\left(-\frac{x^2}{9a}\right), \quad (5.1)$$

$$\mathbb{P}\left(\sup_{0 \leq s \leq t} \Xi(s) < 1\right) \leq 2\exp\left(-\frac{t}{50}\right). \quad (5.2)$$

PROOF. By definition of $\Xi$ in (4.3) (with $\kappa = 0$),

$$\Xi(t) = \int_0^t \sqrt{1 - e^{-\Xi(s)}}\, d\beta(s) + \tfrac{1}{2}\int_0^t e^{-\Xi(s)}\, ds.$$

It follows from the Dubins–Schwarz theorem ([19], Theorem V.1.6) that

$$(5.3) \qquad \Xi(t) = \gamma\left(\int_0^t (1 - e^{-\Xi(s)})\, ds\right) + \tfrac{1}{2}\int_0^t e^{-\Xi(s)}\, ds, \qquad t \geq 0,$$

where $\gamma(\cdot)$ denotes a one-dimensional Brownian motion. Since $\Xi(s) \geq 0$ for all $s \geq 0$, we have

$$\sup_{0 \leq t_1 \leq t_2 < t, t_2 - t_1 < a} |\Xi(t_2) - \Xi(t_1)| \leq \sup_{0 \leq s_1 < s_2 < t, s_2 - s_1 < a} |\gamma(s_2) - \gamma(s_1)| + \frac{a}{2}.$$

According to Lemma 1.1.1 of [5],

$$\mathbb{P}\left(\sup_{0 \leq s_1 < s_2 < t, s_2 - s_1 < a} |\gamma(s_2) - \gamma(s_1)| > \frac{x}{2}\right) \leq c\frac{t}{a}\exp\left(-\frac{x^2}{9a}\right).$$

This implies (5.1).

To prove (5.2), we note that on $\{\sup_{0 \leq s \leq t} \Xi(s) < 1\}$, we have $\int_0^t e^{-\Xi(s)}\, ds \geq t/e$, hence, for all $t > 100$,

$$\left\{\sup_{0 \leq s \leq t} \Xi(s) < 1\right\} \subset \left\{\gamma\left(\int_0^t (1 - e^{-\Xi(s)})\, ds\right) < 1 - \frac{t}{2e} < -\frac{t}{5}\right\}$$

$$\subset \left\{\inf_{0 \leq u \leq t} \gamma(u) < -\frac{t}{5}\right\}.$$

The estimate (5.2) now follows from the usual estimate for Brownian tails. □

LEMMA 5.2. *For any $\varepsilon \in (0,1)$, there exists some $v_0 = v_0(\varepsilon) > 0$ such that for all $x, v > v_0$, we have*

$$\frac{2}{\pi}\exp\left(-(1+\varepsilon)\frac{\pi^2}{8}\frac{v}{x^2}\right) \leq \mathbb{P}\left(\sup_{0 \leq s \leq v} \Xi(s) < x\right) \leq 9\exp\left(-(1-\varepsilon)\frac{\pi^2}{8}\frac{v}{x^2}\right).$$

PROOF. Assume for the moment that $\Xi$ starts from $\Xi(0) = 1$. Let

$$f(x) \stackrel{\text{def}}{=} \int_1^x \frac{dy}{1 - e^{-y}}, \qquad x > 0,$$

be the scale function of $\Xi$. Since $t \mapsto f(\Xi(t))$ is a continuous local martingale, it follows from the Dubins–Schwarz theorem ([19], Theorem V.1.6) that

$$f(\Xi(t)) = B\left(\int_0^t \frac{ds}{1 - e^{-\Xi(s)}}\right), \qquad t \geq 0,$$

for some one-dimensional Brownian motion $B$ starting from 0. Therefore, by writing $T_\Xi(x) = \inf\{s > 0 : \Xi(s) > x\}$, we have

$$T_\Xi(x) = \int_0^{\sigma(f(x))} (1 - e^{-f^{-1}(B(s))})\, ds = \int_{-\infty}^{f(x)} (1 - e^{-f^{-1}(y)}) L(\sigma(f(x)), y)\, dy,$$

where $f^{-1}$ is the inverse of the increasing function $f$, $\sigma(x) \stackrel{\text{def}}{=} \inf\{s : B(s) = x\}$ for $x \in \mathbb{R}$, and $L$ is the local time of $B$.

Observe that $f^{-1}(y) \sim y$ as $y \to \infty$, and $f^{-1}(y) \sim e^{-|y|}$ as $y \to -\infty$. Let $y_0 = y_0(\varepsilon) > 0$ be sufficiently large such that $e^{-f^{-1}(y)} < \varepsilon/2$ for all $y \geq y_0$. Denote by $b(\varepsilon) \stackrel{\text{def}}{=} \sup_{-\infty < y \leq y_0}(1 - e^{-f^{-1}(y)})e^{|y|} < \infty$. Then for all large $x$,

$$(5.4) \quad T_\Xi(x) \geq \left(1 - \frac{\varepsilon}{2}\right) \int_{y_0}^{f(x)} L(\sigma(f(x)), y)\, dy,$$

$$(5.5) \quad T_\Xi(x) \leq \int_{y_0}^{f(x)} L(\sigma(f(x)), y)\, dy + b(\varepsilon) \int_{-\infty}^{y_0} L(\sigma(f(x)), y) e^{-|y|}\, dy.$$

For the lower bound in Lemma 5.2, we note that by (5.4) and (2.3),

$$\mathbb{P}\left\{\sup_{0 \leq s \leq v} \Xi(s) < x\right\} = \mathbb{P}\{T_\Xi(x) > v\}$$

$$\geq \mathbb{P}\left\{\int_{y_0}^{f(x)} L(\sigma(f(x)), y)\, dy > \frac{v}{1 - \varepsilon/2}\right\}$$

$$\geq \frac{2}{\pi} \exp\left(-\frac{\pi^2}{8} \frac{v}{(1 - \varepsilon/2)(f(x) - y_0)^2}\right).$$

Since $f(x) \sim x$, $x \to \infty$, this yields the lower bound in Lemma 5.2 in the case when $\Xi$ starts from $\Xi(0) = 1$, and a fortiori, in the case when $\Xi$ starts from $\Xi(0) = 0$ by a comparison theorem for diffusion processes ([19], Theorem IX.3.7).

For the upper bound in Lemma 5.2, we again assume $\Xi(0) = 1$ for the moment. By (5.5) and the triangular inequality, for $r \geq v_0$,

$$\mathbb{P}\left\{\sup_{0 \leq s \leq r} \Xi(s) < x\right\} = \mathbb{P}\{T_\Xi(x) > r\}$$

$$\leq \mathbb{P}\left\{\int_{y_0}^{f(x)} L(\sigma(f(x)), y)\, dy \geq (1 - \varepsilon/2)r\right\}$$

$$+ \mathbb{P}\left\{\int_{-\infty}^{y_0} L(\sigma(f(x)), y) e^{-|y|}\, dy > \frac{\varepsilon}{2b(\varepsilon)} r\right\}.$$

The first probability expression on the right-hand side is $\leq \frac{4}{\pi} \exp(-\frac{\pi^2(1-\varepsilon/2)r}{8(f(x)-y_0)^2})$ [see (2.3)], whereas the second is $\leq 3\exp(-\frac{\varepsilon}{16b(\varepsilon)f(x)}r) + 2\exp(-\frac{\varepsilon}{8y_0 b(\varepsilon)f(x)}r)$ in light of (2.6). Therefore, if $\Xi(0) = 1$, then for all $r, x \geq v_0$,

$$(5.6) \qquad \mathbb{P}\left\{\sup_{0 \leq s \leq r} \Xi(s) < x\right\} \leq \left(\frac{4}{\pi} + 5\right) \exp\left(-\frac{\pi^2(1-\varepsilon)r}{8x^2}\right).$$

We are now back to the case $\Xi(0) = 0$ we were studying. By (5.2), for any $v > 100/\varepsilon$, $\mathbb{P}\{T_\Xi(1) > \varepsilon v\} \leq 2\exp(-\frac{\varepsilon v}{50})$, which, in view of (5.6) [taking $r \stackrel{\text{def}}{=} (1-\varepsilon)v$ there], yields that

$$\mathbb{P}\left\{\sup_{0 \leq s \leq v} \Xi(s) < x\right\} \leq \left(\frac{4}{\pi} + 5\right) \exp\left(-\frac{\pi^2(1-\varepsilon)^2 v}{8x^2}\right) + 2\exp\left(-\frac{\varepsilon v}{50}\right),$$

which yields the upper bound in Lemma 5.2 [since $\varepsilon \in (0,1)$ is arbitrary]. □

We are now ready to give the proof of the upper bound in Theorem 1.2.

PROOF OF THEOREM 1.2: THE UPPER BOUND. Observe that, by (4.2),

$$\mathbb{P}\left\{\sup_{0 \leq s \leq t} X(s) > v\right\} = \mathbb{P}(H(v) < t) \leq \mathbb{P}\{\Theta_1(v) < t\}.$$

By Theorem 4.1, $\Theta_1(v)$ is distributed as $4\int_0^v e^{\Xi(s)}\,ds - 4v$. Thus,

$$(5.7) \qquad \mathbb{P}\left\{\sup_{0 \leq s \leq t} X(s) > v\right\} \leq \mathbb{P}\left\{\int_0^v e^{\Xi(s)}\,ds < \frac{t}{4} + v\right\}.$$

According to (5.1) (taking $x \stackrel{\text{def}}{=} 3$ and $a \stackrel{\text{def}}{=} \frac{1}{v}$ there), we have

$$\mathbb{P}\left\{\sup_{0 \leq t_1 \leq t_2 \leq v, t_2 - t_1 < 1/v} |\Xi(t_2) - \Xi(t_1)| > 3\right\} \leq cv^2 e^{-v},$$

where $c$ is the numerical constant in (5.1). On the event $\{\sup_{0 \leq t_1 \leq t_2 \leq v, t_2-t_1 < 1/v} |\Xi(t_2) - \Xi(t_1)| \leq 3\}$, we have $\int_0^v e^{\Xi(s)}\,ds \geq \frac{1}{v}\exp(\sup_{0 \leq s \leq v} \Xi(s) - 3)$. Plugging this into (5.7) yields that for all sufficiently large $v$ and $t$,

$$(5.8) \quad \begin{aligned} \mathbb{P}\left\{\sup_{0 \leq s \leq t} X(s) > v\right\} &\leq cv^2 e^{-v} + \mathbb{P}\left\{\sup_{0 \leq s \leq v} \Xi(s) - 3 < \log\left(\frac{tv}{4} + v^2\right)\right\} \\ &\leq cv^2 e^{-v} + 9\exp\left(-(1-2\varepsilon)\frac{\pi^2}{8}\frac{v}{\log^2(tv)}\right), \end{aligned}$$

the last inequality being a consequence of the upper bound in Lemma 5.2. We mention that (5.8) was already used in Section 3 to prove the estimate (3.18).

Since $\log v = o(\log t)$, (5.8) yields the upper bound in Theorem 1.2. □

5.2. *Theorem* 1.2: *the lower bound.* The ideas in this section essentially go back to [2]. For $a, b \in \mathbb{R}$, we define

$$\overline{W}(a,b) \stackrel{\text{def}}{=} \sup_{0 \leq s \leq 1} W(a + s(b-a)),$$

$$\underline{W}(a,b) \stackrel{\text{def}}{=} \inf_{0 \leq s \leq 1} W(a + s(b-a)),$$

$$W^{\#}(a,b) \stackrel{\text{def}}{=} \sup_{0 \leq s \leq t \leq 1} [W(a + t(b-a)) - W(a + s(b-a))].$$

Note that $W^{\#}(a,b) \neq W^{\#}(b,a)$ in general. Let $\mathbf{P}_{x,\omega}$ be the quenched probability under which the diffusion $X$ starts from $x$.

Recall that $H(y) = \inf\{t \geq 0 : X(t) = y\}$. Let

$$(5.9) \qquad \psi(x) := \mathbb{P}\left\{\inf_{|y| \leq 1/2} L(\sigma(1) \wedge \sigma(-1), y) < x\right\}, \qquad x > 0.$$

We start with the following lemma. We mention that $(W(y), a \leq y \leq c)$ is *not* necessarily a valley in the sense of Brox [2].

LEMMA 5.3. *Let $a < x < c$ and let $\lambda > 0$.*

1. *We have*

$$(5.10) \quad \mathbf{P}_{x,\omega}\{H(a) \wedge H(c) > \lambda(c-a)^2 e^{\min(W^{\#}(a,c), W^{\#}(c,a))}\} \leq \frac{24}{\lambda} + \frac{96}{\lambda^2}.$$

2. *Let*

$$(5.11) \qquad\qquad a' = \sup\{y \leq x : W(y) - W(x) > \lambda\},$$

$$(5.12) \qquad\qquad c' = \inf\{y > x : W(y) - W(x) > \lambda\}.$$

*If*

$$(5.13) \ (c'-a')e^{\lambda} \leq \tfrac{1}{2}\min\left\{\int_a^x e^{W(y)-W(x)}\,dy, \int_x^c e^{W(y)-W(x)}\,dy\right\} \stackrel{\text{def}}{=} \tfrac{1}{2}\Gamma(a,c),$$

*then for all $0 < \varepsilon < 1$,*

$$(5.14) \qquad \mathbf{P}_{x,\omega}\{H(a) \wedge H(c) \leq \varepsilon(c'-a')e^{-\lambda}\Gamma(a,c)\} \leq \psi(\varepsilon),$$

*where $\psi(\cdot)$ is defined in (5.9).*

PROOF. 1. According to Brox ([2], pages 1213 and 1214, proof of (i); we mention that the $m$ in [2] is our $c$, and the $\alpha$ in [2] is 1 here),

$$H(a) \wedge H(c) \leq (c-a)^2 e^{W^{\#}(a,c)} \Delta_4,$$

where $\Delta_4 \stackrel{\text{law}}{=} \sup_{-\infty < y \leq 1} L(\sigma(1), y) + \Delta_5$ and $\mathbf{E}_{x,\omega}(\Delta_5^2) \leq 12$. A similar estimate holds with $W^{\#}(c,a)$ instead of $W^{\#}(a,c)$.

Therefore, (5.10) will follow from Chebyshev's inequality once we can show that for all $\lambda > 0$,

$$(5.15) \qquad \mathbb{P}\left\{\sup_{-\infty < y \leq 1} L(\sigma(1), y) > \lambda\right\} \leq \frac{6}{\lambda}.$$

To prove (5.15), we note that by Fact 2.1, $y \in [0,1] \mapsto L(\sigma(1), y)$ is a BESQ(2) starting from 0, and $y \in [0, \infty) \mapsto L(\sigma(1), -y)$ is a BESQ(0) starting from $L(\sigma(1), 0)$. Using successively the triangular inequality, the reflection principle for BESQ(2) and the martingale property of BESQ(0), we obtain that

$$\mathbb{P}\left\{\sup_{-\infty < y \leq 1} L(\sigma(1), y) > \lambda\right\}$$

$$\leq \mathbb{P}\left\{\sup_{0 \leq y \leq 1} L(\sigma(1), y) > \lambda\right\} + \mathbb{E}(\mathbf{1}_{\{L(\sigma(1),0) < \lambda\}} \mathbf{1}_{\{\sup_{-\infty < y \leq 0} L(\sigma(1), y) > \lambda\}})$$

$$\leq 2\mathbb{P}\{L(\sigma(1), 0) > \lambda\} + \mathbb{E}\left(\mathbf{1}_{\{L(\sigma(1),0) < \lambda\}} \frac{L(\sigma(1), 0)}{\lambda}\right).$$

Since $L(\sigma(1), 0)$ has the exponential distribution of mean 2, this leads to

$$\mathbb{P}\left\{\sup_{-\infty < y \leq 1} L(\sigma(1), y) > \lambda\right\} \leq 2e^{-\lambda/2} + \frac{2}{\lambda}\int_0^{\lambda/2} y e^{-y}\,dy \leq \frac{6}{\lambda},$$

proving (5.15) and, thus, (5.10).

2. The proof of (5.14) is essentially from [2], page 1215, line -11. Without loss of generality, we assume $x = 0$. Under (5.13), we have $a < a' < 0 < c' < c$. In view of (4.1), we have, by the occupation time formula,

$$H(a) \wedge H(c) = T(\sigma(A(a)) \wedge \sigma(A(c)))$$

$$= \int_a^c e^{-W(y)} L(\sigma(A(a)) \wedge \sigma(A(c)), A(y))\,dy$$

$$\geq (c' - a')e^{-\lambda} \inf_{a' \leq y \leq c'} L(\sigma(A(a)) \wedge \sigma(A(c)), A(y))$$

$$\geq (c' - a')e^{-\lambda} \inf_{|z| \leq (c'-a')e^\lambda} L(\sigma(A(a)) \wedge \sigma(A(c)), z),$$

since $A(c') \leq c'e^\lambda$ and $A(a') \geq -|a'|e^\lambda$. In view of (5.13) and scaling,

$$\mathbf{P}_{x,\omega}\{H(a) \wedge H(c) \leq \varepsilon(c'-a')e^{-\lambda}\Gamma(a,c)\}$$

$$\leq \mathbb{P}\left\{\inf_{|y| \leq 1/2} L(\sigma(1) \wedge \sigma(-1), y) < \varepsilon\right\},$$

proving (5.14). $\square$

We now have all the ingredients to prove the lower bound in Theorem 1.2.

PROOF OF THEOREM 1.2: THE LOWER BOUND. Fix a small $\varepsilon > 0$. Consider large $t$ and $v$ such that $t^\varepsilon \geq v \geq \log^2 t$. Let $r = \frac{\log t}{1-10\varepsilon}$. For $r > 0$, define $d_-(r) \stackrel{\text{def}}{=} \sup\{t < 0 : W^\#(t, 0) > r\}$.

Define three random times $\alpha > m > \eta > v$ by

$$\eta = \inf\{s > v : W(s) - W(v) = -(1-3\varepsilon)r\},$$
$$\alpha = \inf\{s > \eta : W(s) - W(\eta) = r\},$$
$$m = \inf\{s > \eta : W(s) = \underline{W}(v, \alpha)\}.$$

We consider the following events concerning the Brownian potential $W$:

$$F_1 \stackrel{\text{def}}{=} \left\{|d_-(r)| < r^2; |\underline{W}(d_-(r), 0)| \leq \varepsilon r; \int_{d_-(r)}^0 e^{W(z)}\,dz > e^{r/2}\right\},$$

$$F_2 \stackrel{\text{def}}{=} \left\{W^\#(0, v) < (1-20\varepsilon)r; \overline{W}(0, v) < \frac{r}{3}\right\},$$

$$F_3 \stackrel{\text{def}}{=} \left\{\eta - v \leq r^{5/2}; \sup_{v \leq s \leq \eta}(W(s) - W(v)) < \varepsilon r; W^\#(v, \eta) < \frac{r}{3};\right.$$
$$\left.\int_v^\eta e^{W(x) - W(\eta)}\,dx > e^{(1-4\varepsilon)r}\right\},$$

$$F_4 \stackrel{\text{def}}{=} \left\{\alpha - \eta \leq r^{5/2}; \inf_{\eta \leq s \leq \alpha}(W(s) - W(\eta)) > -\varepsilon r;\right.$$
$$W^\#(\eta, m) < \frac{r}{3}; W^\#(\alpha, \eta) < \frac{r}{3};$$
$$\left.\sup_{m \leq x \leq m+1}(W(x) - W(m)) < r^{2/3}; \int_m^\alpha e^{W(x) - W(m)}\,dx > e^{(1-\varepsilon)r}\right\}.$$

Observe that by the strong Markov property, the events $(F_j)_{1 \leq j \leq 4}$ are independent. Moreover, $P\{F_3\}$ and $P\{F_4\}$ do not depend on $v$.

Clearly,

$$\mathbf{P}_{0,\omega}\{X(t) > v\} \geq \mathbf{P}_{0,\omega}\{H(m) < t; X(t) > v\}$$
(5.16)
$$\geq \mathbf{P}_{0,\omega}\{H(m) < t\}\mathbf{P}_{m,\omega}\{H(v) \wedge H(\alpha) > t\}.$$

We have

$$\mathbf{P}_{0,\omega}\{H(m) \geq t\} \leq \mathbf{P}_{0,\omega}\{H(d_-(r)) < H(m)\} + \mathbf{P}_{0,\omega}\{H(d_-(r)) \wedge H(m) \geq t\}$$
$$= \frac{A(m)}{A(m) - A(d_-(r))} + \mathbf{P}_{0,\omega}\{H(d_-(r)) \wedge H(m) \geq t\}.$$

Let $\omega \in \bigcap_{j=1}^{4} F_j$. Since $A(m) \leq m e^{\overline{W}(m)} \leq m e^{r/3+\varepsilon r}$, $|A(d_-(r))| \geq e^{r/2}$, and $W^{\#}(d_-(r), m) \leq (1-19\varepsilon)r$, we apply (5.10) to $(d_-(r), 0, m)$ and arrive at

$$\mathbf{P}_{0,\omega}\{H(m) \geq t\} \leq \frac{me^{r/3+\varepsilon r}}{e^{r/2}} + \frac{24}{t}(m - d_-(r))^2 e^{(1-19\varepsilon)r}$$
$$+ \frac{96}{t^2}(m - d_-(r))^4 e^{2(1-19\varepsilon)r}.$$

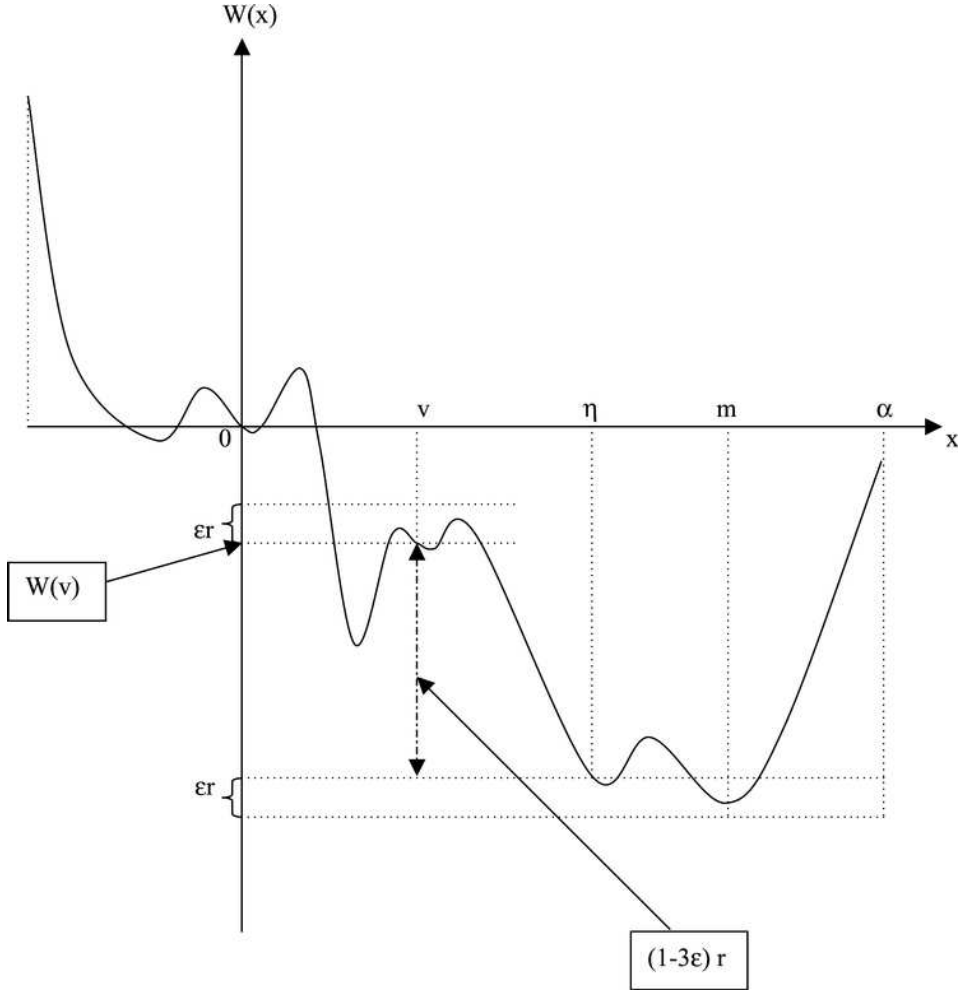

Fig. 1.

Note that $m \leq \alpha \leq r^3 + v \leq 2t^\varepsilon$ and $t = e^{(1-10\varepsilon)r}$. We obtain that

$$\mathbf{P}_{0,\omega}\{H(m) \geq t\} \leq e^{-6\varepsilon r} \qquad \text{if } \omega \in \bigcap_{j=1}^{4} F_j. \tag{5.17}$$

To apply (5.14) to $(v, m, \alpha)$, we choose $\lambda = r^{2/3}$ and verify that the assumption (5.13) is satisfied on $\bigcap_{j=1}^{4} F_j$, because

$$\Gamma(v, \alpha) = \min\left\{\int_v^m e^{W(y)-W(m)}\, dy, \int_m^\alpha e^{W(y)-W(m)}\, dy\right\}$$

$$\geq \min\left\{\int_v^\eta e^{W(y)-W(\eta)}\, dy, e^{(1-\varepsilon)r}\right\}$$

$$\geq e^{(1-4\varepsilon)r} \geq 2(c'-a')e^\lambda.$$

It follows from (5.14) that

$$\mathbf{P}_{m,\omega}\{H(v) \wedge H(\alpha) \leq t\} \leq \psi\left(\frac{te^\lambda}{(c'-a')\Gamma(v,\alpha)}\right).$$

Since $c' - a' \geq c' - m \geq 1$ by definition of $F_4$,

$$\psi\left(\frac{te^\lambda}{(c'-a')\Gamma(v,\alpha)}\right) \leq \psi(e^{(1-10\varepsilon)r + r^{2/3} - (1-4\varepsilon)r}) \leq \psi(e^{-5\varepsilon r}).$$

Therefore, for all large $r \geq r_0(\varepsilon)$, we get from (5.16) that

$$\mathbf{P}_{0,\omega}\{X(t) > v\} \geq (1 - e^{-6\varepsilon r})(1 - \psi(e^{-5\varepsilon r})) \geq \tfrac{1}{2} \qquad \text{if } \omega \in \bigcap_{j=1}^{4} F_j.$$

Hence,

$$P\{X(t) > v\} = E(\mathbf{P}_{0,\omega}\{X(t) > v\}) \geq \tfrac{1}{2} P\left\{\bigcap_{j=1}^{4} F_j\right\} = \tfrac{1}{2} \prod_{j=1}^{4} P\{F_j\}, \tag{5.18}$$

by the independence of $F_j$. When $r \to \infty$, $P\{\eta - v > r^{5/2}\} \to 0$ and

$$\frac{1}{r} \log \int_v^\eta e^{W(x)-W(\eta)}\, dx \sim \frac{1}{r} \sup_{v \leq x \leq \eta}(W(x) - W(\eta)) \geq (1 - 3\varepsilon).$$

It follows that

$$\liminf_{r \to \infty} P\{F_3\} \geq \liminf_{r \to \infty} P\left\{\sup_{v \leq s \leq \eta}(W(s) - W(v)) < \varepsilon r; W^\#(v,\eta) < \frac{r}{3}\right\}$$

$$= C(\varepsilon) > 0,$$

for some constant $C = C(\varepsilon)$ depending only on $\varepsilon$. The same holds for $P\{F_4\}$ and $P\{F_1\}$. It follows that for all large $r \geq r_0$, we have

$$(5.19) \qquad P\{F_1\}P\{F_3\}P\{F_4\} \geq C'(\varepsilon) > 0.$$

Finally, we recall the asymptotic expansion of the distribution of $(W^\#(0, v), \overline{W}(0, v))$ ([9], Theorem 2.1): for any fixed $0 < a \leq 1$, when $\delta \to 0+$,

$$P\{W^\#(0,1) < \delta; \overline{W}(0,1) < a\delta\} \sim \frac{4\sin(a\pi/2)}{\pi} \exp\left(-\frac{\pi^2}{8\delta^2}\right).$$

It follows from scaling that when $r \to \infty$,

$$P\{F_2\} \sim \frac{4\sin(\pi/6(1-20\varepsilon))}{\pi} \exp\left(-\frac{\pi^2}{8(1-20\varepsilon)^2}\frac{v}{r^2}\right).$$

Plugging this into (5.18) and (5.19) implies

$$\liminf_{t,v\to\infty, v \gg \log^2 t, \log v = o(\log t)} \frac{\log^2 t}{v} \log \mathbb{P}\{X(t) > v\} \geq -\frac{\pi^2}{8}.$$

The lower bound in Theorem 1.2 is proved. $\square$

## REFERENCES

[1] BASS, R. F. and GRIFFIN, P. S. (1985). The most visited site of Brownian motion and simple random walk. *Z. Wahrsch. Verw. Gebiete* **70** 417–436. MR803682
[2] BROX, T. (1986). A one-dimensional diffusion process in a Wiener medium. *Ann. Probab.* **14** 1206–1218. MR866343
[3] COMETS, F. and POPOV, S. (2003). Limit law for transition probabilities and moderate deviations for Sinai's random walk in random environment. *Probab. Theory Related Fields* **126** 571–609. MR2001198
[4] CSÁKI, E. and FÖLDES, A. (1987). A note on the stability of the local time of a Wiener process. *Stochastic Process. Appl.* **25** 203–213. MR915134
[5] CSÖRGŐ, M. and RÉVÉSZ, P. (1981). *Strong Approximations in Probability and Statistics*. Academic Press, New York. MR666546
[6] FELLER, W. (1971). *An Introduction to Probability Theory and Its Applications*, **2**, 2nd ed. Wiley, New York.
[7] GANTERT, N. and ZEITOUNI, O. (1999). Large deviations for one-dimensional random walk in a random environment—a survey. In *Random Walks* (P. Révész and B. Tóth, eds.) 127–165. Bolyai Math. Soc., Budapest. MR1752893
[8] GETOOR, R. K. and SHARPE, M. J. (1979). Excursions of Brownian motion and Bessel processes. *Z. Wahrsch. Verw. Gebiete* **47** 83–106. MR521534
[9] HU, Y. and SHI, Z. (1998). The limits of Sinai's simple random walk in random environment. *Ann. Probab.* **26** 1477–1521. MR1675031
[10] HU, Y., SHI, Z. and YOR, M. (1999). Rates of convergence of diffusions with drifted Brownian potentials. *Trans. Amer. Math. Soc.* **351** 3915–3934. MR1637078
[11] KARLIN, S. and TAYLOR, H. M. (1981). *A Second Course in Stochastic Processes*. Academic Press, New York. MR611513
[12] KAWAZU, K. and TANAKA, H. (1997). A diffusion process in a Brownian environment with drift. *J. Math. Soc. Japan* **49** 189–211. MR1601361

[13] Kent, J. (1978). Some probabilistic properties of Bessel functions. *Ann. Probab.* **6** 760–770. MR501378
[14] Kesten, H. (1965). An iterated logarithm law for local time. *Duke Math. J.* **32** 447–456. MR178494
[15] Lamperti, J. (1972). Semi-stable Markov processes, I. *Z. Wahrsch. Verw. Gebiete* **22** 205–225. MR307358
[16] Le Doussal, P., Monthus, C. and Fisher, D. S. (1999). Random walkers in one-dimensional random environments: Exact renormalization group analysis. *Phys. Rev. E* **59** 4795–4840. MR1682204
[17] Pitman, J. and Yor, M. (1982). A decomposition of Bessel bridges. *Z. Wahrsch. Verw. Gebiete* **59** 425–457. MR656509
[18] Révész, P. (1990). *Random Walk in Random and Non-Random Environments.* World Scientific, Singapore. MR1082348
[19] Revuz, D. and Yor, M. (1999). *Continuous Martingales and Brownian Motion*, 3rd ed. Springer, Berlin. MR1725357
[20] Schumacher, S. (1985). Diffusions with random coefficients. *Contemp. Math.* **41** 351–356. MR814724
[21] Shi, Z. (2001). Sinai's walk via stochastic calculus. In *Panoramas et Synthèses* **12** (F. Comets and E. Pardoux, eds.). Société Mathématique de France.
[22] Sinai, Ya. G. (1982). The limiting behaviours of a one-dimensional random walk in a random medium. *Theory Probab. Appl.* **27** 256–268. MR657919
[23] Taleb, M. (2001). Large deviations for a Brownian motion in a drifted Brownian potential. *Ann. Probab.* **29** 1173–1204. MR1872741
[24] Tanaka, H. (1995). Diffusion processes in random environments. In *Proc. of the International Congress of Mathematicians* **2** (S. D. Chatterji, ed.) 1047–1054. Birkhäuser, Basel. MR1404004
[25] Warren, J. and Yor, M. (1998). The Brownian burglar: Conditioning Brownian motion by its local time process. *Séminaire de Probabilités XXXII. Lecture Notes in Math.* **1686** 328–342. Springer, Berlin. MR1655303
[26] Zeitouni, O. (2001). Lecture notes on random walks in random environment. *Saint Flour Lecture Notes*.

Laboratoire de Probabilités
et Modèles Aléatoires
CNRS UMR 7599
Université Paris VI
4 Place Jussieu
F-75252 Paris Cedex 05
France
e-mail: hu@ccr.jussieu.fr
e-mail: zhan@proba.jussieu.fr